\documentclass[notitlepage,twoside,a4paper]{amsart}
\usepackage{mathrsfs}
\usepackage{amsmath,amssymb,enumerate}
\usepackage{epsfig,fancyhdr,color}
\usepackage{epstopdf}
\usepackage{amssymb}
\usepackage{amsmath,amsthm}
\usepackage{latexsym}
\usepackage{amscd}
\usepackage{psfrag}
\usepackage{graphicx}
\usepackage{epsf}
\usepackage[latin1]{inputenc}
\usepackage[all]{xy}
\usepackage{tikz}
\usepackage{prettyref}
\usepackage{subfigure}
\usepackage{float}
\usepackage{color}
\usepackage{array}
\usepackage{hyperref}
\usepackage{cite}
\usepackage[T1]{fontenc}

\usepackage[totalwidth=16cm,totalheight=24cm]{geometry}

\newcommand{\II}{I\hspace{-0.1cm}I}
\newcommand{\III}{I\hspace{-0.1cm}I\hspace{-0.1cm}I}

\newtheorem{theorem}{\rm\bf Theorem}[section]
\newtheorem*{theorem*}{\rm\bf Theorem A}
\newtheorem*{theorem*'}{\rm\bf Theorem A*}
\newtheorem*{theorem**}{\rm\bf Theorem B}
\newtheorem*{theorem***}{\rm\bf Theorem C}
\newtheorem*{theorem****}{\rm\bf Theorem D}
\newtheorem*{question*}{\rm\bf Question}
\newtheorem*{remark*}{\rm\bf Remark}

\newtheorem{proposition}[theorem]{\rm\bf Proposition}

\newtheorem{corollary}[theorem]{\rm\bf Corollary}
\newtheorem{definition}[theorem]{\rm\bf Definition}
\newtheorem{remark}[theorem]{\rm\bf Remark}

\newtheorem{question}[theorem]{\rm\bf Question}

\newtheorem{lemma}[theorem]{\rm\bf Lemma}

\newtheoremstyle{named}{}{}{\itshape}{}{\bfseries}{.}{.5em}{#1 \thmnote{#3}}
\theoremstyle{named}

\newcommand{\PSL}{\operatorname{PSL}}

\newcommand{\CP}{{\mathbb CP}}

\newcommand{\RP}{{\mathbb {RP}}}

\newcommand{\AdS}{\mathbb{ADS}}

\newcommand{\op}{\Omega^+_{\mathcal{C}}}

\newcommand{\CH}{\rm{CH}}

\newcounter{notes}%

\def\interieur#1{\mathord{\mathop{\kern 0pt #1}\limits^\circ}}

\title[]{Embedded Convex Surfaces in Hyperbolic and Anti-de Sitter Spaces}
\author{Abderrahim Mesbah}
\address{Abderrahim Mesbah \newline
Beijing Institute of Mathematical Sciences and Applications, Beijing, China\\
}
\email{abderrahimmesbah@bimsa.cn}

\date{v0, \today}

\begin{document}
\maketitle
\medskip
\begin{abstract}
We show that given a quasi-circle $C$ in $\partial_{\infty}\mathbb{H}^3$ (respectively in $\partial_{\infty}\AdS^3$) and a complete conformal metric $h$ on $\mathbb{D}$ whose curvature $K_h$ takes values in a compact subset of $(-1,0)$ (respectively $(-\infty,-1)$), with all derivatives bounded with respect to the hyperbolic metric, there exists a smooth isometric embedding $V : (\mathbb{D}, h) \to \mathbb{H}^3$ (respectively $V : (\mathbb{D}, h) \to \AdS^3$) such that $V$ extends continuously to a homeomorphism $\partial V : S^1 \to C$. In the case of hyperbolic space, the statement still holds if $C$ is a Jordan curve.
\end{abstract}
\section{Introduction}
A classical problem in differential geometry, known as the \emph{Weyl problem}, was proposed by Weyl in 1915. The problem asks which metrics on the sphere can be realized as the induced metric on the boundary of a compact convex subset of $\mathbb{R}^3$.\\
The first major progress on this problem was made by Lewy \cite{Lewy}, who established the existence of such a realization in the analytic category, while the corresponding uniqueness result in the analytic category is due to Cohn-Vossen \cite{Cohn-Vossen}.\\
A complete solution to the Weyl problem was later obtained by Nirenberg, who proved existence, and Herglotz, who proved uniqueness. On the other hand, this solution can also be deduced from Alexandrov's realization theorem for intrinsic metrics of convex surfaces \cite{alexandrov-weyl}, together with a deep result of Pogorelov \cite{Pogorelov1973} showing that smooth metrics admit only smooth embeddings. These results yield the following theorem.
\begin{theorem}[Weyl problem]\label{intoduction-thoerem}
Any smooth metric of positive Gaussian curvature on the sphere is induced on the boundary of a unique (up to Euclidean isometries) smooth strictly convex body in $\mathbb{R}^3$.
\end{theorem}

Moreover, the work of Alexandrov \cite{alexandrov-weyl} and Pogorelov \cite{Pogorelov1973} extended Theorem~\ref{intoduction-thoerem} to hyperbolic space.

\begin{theorem}[Alexandrov, Pogorelov]
Any smooth metric on the sphere with curvature $K > -1$ is induced on the boundary of a unique (up to isometry) convex subset of $\mathbb{H}^3$ with smooth boundary.
\end{theorem}
Another setting in which a similar question can be asked is that of hyperbolic manifolds. In \cite{Labourie-prescribed}, Labourie characterized the metrics that can be realized on the boundary of a compact hyperbolic manifold. Later, Schlenker \cite{Schlenker2006} proved the rigidity of such manifolds realizing a given metric.\\ 
Labourie \cite{Labourie-curvature} also showed that the complement of the convex core of a convex co-compact hyperbolic manifold (or, more precisely, a hyperbolic end, see \cite[Definition 1.2]{Labourie-curvature}) is foliated by surfaces of constant Gaussian curvature, called $K$-surfaces, each leaf has constant curvature equal to $K$, and $K$ varies in $(-1,0)$. Moreover, he proved that given a smooth function $K : S \to (-1,0)$ and a hyperbolic end $E = S \times (-1,0)$ (where $S$ is supposed to be a closed surface of genus bigger or equal to $2$), there exists a unique immersion of $S$ whose image is an incompressible surface with induced metric of curvature $K$. Later, Rosenberg and Spruck in \cite{R-S} obtained an analogous result in the universal cover $\mathbb{H}^3$. They showed that for any Jordan curve $C \subset \partial_{\infty}\mathbb{H}^3$, the complement $\mathbb{H}^3 \setminus CH(C)$, where $CH(C)$ denotes the convex hull of $C$, is foliated by $K$-surfaces. Each leaf of this foliation is a topological disk that spans $C$ at the boundary at infinity and has constant Gaussian curvature equal to $K$, where $K$ varies in $(-1,0)$.\\ 
Similar results to those in \cite{Labourie-prescribed}, concerning the universal cover, can be found in \cite{bonsante2021induced}. The authors show that for any quasisymmetric homeomorphism $f$, and for any $K\in(-1,0)$, there exist a quasicircle $C$ and two convex embedded discs $S^{\pm}$ with asymptotic boundary $\partial_{\infty}S^{\pm}=C$, such that the induced metric on each $S^{\pm}$ is isometric to $\mathbb{D}$ endowed with the complete conformal metric of constant curvature $K$, and the corresponding gluing map is equal to $f$. Later, in \cite{Chen-Schlenker}, the authors extended this result to a larger class of negatively curved metrics. Unlike the group-invariant setting, no uniqueness statement is currently known in the case of the universal cover.
\\
These results naturally lead to existence questions for isometric embeddings with a prescribed conformal metric on the unit disk $\mathbb{D}$. Schlenker \cite[Question 4.9]{JM-Weyl} posed the following question:
\begin{question}
Let $h$ be a complete conformal metric on the disk $\mathbb{D}$ with curvature $K > -1$. Is there a unique isometric embedding of $(\mathbb{D}, h)$ into $\mathbb{H}^3$ such that the ideal boundary of the image is a round circle?
\end{question}
In this work, we provide a partial answer to this question. Before stating our main theorem, we would like to emphasize an important difference between our result and those of \cite{bonsante2021induced} and \cite{Chen-Schlenker}. In our setting, we prove the existence of convex surfaces with a prescribed metric and a prescribed asymptotic boundary, for an arbitrary Jordan curve. By contrast, the results of \cite{bonsante2021induced} and \cite{Chen-Schlenker} establish the existence of convex domains spanning a given quasicircle and realizing a prescribed gluing map. In particular, once the gluing map is fixed, the corresponding convex domain (and hence the quasicircle) is expected to be unique. Our theorem does not involve a prescribed gluing map, and therefore should be viewed as addressing a different existence problem.
\begin{theorem}\label{first}
Let $h$ be a complete conformal metric on the disc $\mathbb{D}$, such that the curvature $K_{h}$ of $h$ takes values in a compact subset of $(-1,0)$, and such that all derivatives of $K_{h}$ are bounded with respect to the hyperbolic metric, at any order. Let $C$ be a Jordan curve in $\partial_{\infty}\mathbb{H}^{3}$.\\  
Then there exists an isometric embedding $V: (\mathbb{D},h) \to \mathbb{H}^3$
that extends continuously to $\partial V: S^1 \to C$.\\
\end{theorem}
From the proof of Theorem \ref{first}, one can observe the existence of at least two such surfaces, each lying in a distinct component of $\mathbb{H}^3 \setminus \CH(C)$, where $\CH(C)$ denotes the hyperbolic convex hull of $C$. In \cite{R-S}, Rosenberg and Spruck show that when the curvature of $h$ is constant, there exist exactly two such surfaces, with each component of $\mathbb{H}^3 \setminus \CH(C)$ containing one of them. This naturally leads to the following question:
\begin{question}
In the general case, do there exist exactly two embedded surfaces as in Theorem \ref{first}, one in each component of $\mathbb{H}^3 \setminus \CH(C)$?
\end{question}
A related setting is considered by Sui in \cite{Sui-Zhenan}, where the author proves the existence of a hypersurface in $\mathbb{H}^n$ with a prescribed scalar curvature and a given smooth submanifold of $\partial_{\infty}\mathbb{H}^n$ as its asymptotic boundary. In \cite{Sui-Wei} and \cite{Sui-Wei-Tamburelli}, Sui and Sun establish the existence of a smooth, complete hypersurface in $\mathbb{H}^n$ with a prescribed Weingarten curvature in the weak sense and asymptotic boundary at infinity.\\ 
An analogous problem arises in the anti-de Sitter space. The three-dimensional anti-de Sitter space $\AdS^3$ can be seen as the Lorentzian analogue of hyperbolic space. Its ideal boundary $\partial_{\infty}\AdS^3$ is identified with $S^1 \times S^1$ (see Section \ref{anti-de sitter}). A quasi-circle in $\partial_{\infty}\AdS^3$ is defined as the graph of a quasi-symmetric homeomorphism of $S^1$.\\
The geometry of $\AdS^3$ shares many features with hyperbolic space. In particular, there is the notion of globally hyperbolic maximal compact manifolds (see Section \ref{Globally-hyperbolic}), which can be viewed as the Lorentzian analogue of quasi-Fuchsian manifolds (see Section \ref{quasi-Fuchsian-manifolds}).\\ 
A Lorentzian analogue of \cite{Labourie-prescribed} is shown by Tamburelli \cite{tamburelli}, where Tamburelli characterized the metrics that can be realized on the boundary of an $\AdS^3$ globally hyperbolic maximal compact manifold. In contrast to the hyperbolic setting, the uniqueness problem turns out to be much more challenging. Nevertheless, Prosanov and Schlenker \cite{Prosanov-Schlenker} provided a partial answer to the uniqueness question. They showed that when the holonomy of the $\AdS$ manifold is sufficiently close to a Fuchsian representation, the manifold is uniquely determined by the induced metrics on its boundary, which are prescribed by Tamburelli's theorem.\\
As in the hyperbolic case, Barbot, B\'eguin, and Zeghib in \cite{curvature-zeghib} obtained results analogous to those of Labourie \cite{Labourie-prescribed}. They showed that any anti-de Sitter end which is diffeomorphic to $S \times (-\infty,-1)$, and $S$ is a closed surface of genus bigger or equal to $2$ (see Section \ref{anti-de sitter}), is foliated by $K$-surfaces, where each leaf has a constant Gaussian curvature equal to $K$, and $K$ varies in $(-\infty,-1)$. They also proved that for any smooth function $K : S \to (-\infty,-1)$, there exists an embedding of $S$ into the $\AdS^3$ end whose induced metric has curvature $K$.\\
Later, Bonsante and Seppi \cite{BS-ksurfaces} established a universal version of the theorem by Barbot, B\'eguin, and Zeghib. They proved that for any quasi-circle $C \subset \partial_{\infty}\AdS^3$, the complement $D(C) \setminus CH(C)$ is foliated by $K$-surfaces, where $D(C)$ denotes the domain of dependence of $C$ (see Section \ref{Globally-hyperbolic}), and $K$ varies in $(-\infty,-1)$. The same authors also obtained analogous results in the Minkowski setting \cite{BS-Minkowski}, showing in particular that one can construct surfaces with prescribed curvature.
More precisely, let $\varphi : \partial \mathbb{D} \to \mathbb{R}$ be a lower semi-continuous and bounded function, and let $\psi : \mathbb{D} \to [a,b]$ for some $0 < a < b < +\infty$. Then there exists a unique spacelike graph $S \subset \mathbb{R}^{2,1}$ whose support function $u$ extends $\varphi$ and whose curvature function is $\psi$.\\
Again as in the hyperbolic setting, the authors of \cite{bonsante2021induced} proved a group-action-free analogue of the result of \cite{tamburelli}. More precisely, they showed that for any quasisymmetric homeomorphism $f$, and for any $K\in(-\infty,-1)$, there exist a quasicircle $C\subset\partial_{\infty}\AdS^3$ and two convex embedded discs $S^{\pm}$ with asymptotic boundary $\partial_{\infty}S^{\pm}=C$, such that the induced metric on each $S^{\pm}$ is isometric to $\mathbb{D}$ endowed with the complete conformal metric of constant curvature $K$, and the corresponding gluing map is equal to $f$. Later, in \cite{curvature-abder}, this result was extended to a larger class of negatively curved metrics. As in the hyperbolic setting, no uniqueness statement is currently known in this setting.\\
In this work, we establish an existence result in the $\AdS^3$ setting analogous to Theorem \ref{first}. We show that for any quasi-circle $C \subset \partial_{\infty}\AdS^3$ and any complete conformal metric $h$ on $\mathbb{D}$ whose curvature lies in the interval $(-\frac{1}{\epsilon}, -1-\epsilon)$ and whose hyperbolic derivatives of all orders are bounded, there exists an isometric embedding of $(\mathbb{D}, h)$ into $\AdS^3$ that spans $C$ at the boundary at infinity.
\begin{theorem}\label{secend}
Let $h$ be a complete, conformal metric on the disc $\mathbb{D}$, such that the sec curvature $K_{h}$ of $h$ takes values in a compact subset of $(-\infty,-1)$, and such that all derivatives of $K_{h}$ are bounded with respect to the hyperbolic metric, at any order.  Let $C$ be a quasi-circle in $\partial_{\infty}\AdS^{3}$.\\  
Then there exists an isometric embedding $V: (\mathbb{D},h) \to \AdS^{3}$
that extends continuously to $\partial V: S^1 \to C$
\end{theorem}
Exactly as in the hyperbolic case, one can notice from our proof that there exist at least two such surfaces, each lying in a distinct component of $D(C) \setminus \CH(C)$, where $D(C)$ denotes the domain of dependence of $C$ (see Definition \ref{domain-of-dependance}) and $\CH(C)$ is the anti-de Sitter convex hull of $C$. Bonsante and Seppi in \cite{BS-ksurfaces} showed that when the curvature of $h$ is constant, there exist exactly two such surfaces, each contained in a different component of $D(C) \setminus \CH(C)$. This naturally leads to the following question:
\begin{question}
In the general case, do there exist exactly two embedded surfaces as in Theorem \ref{secend}, one in each component of $D(C) \setminus \CH(C)$?
\end{question}
\section{outline of the paper}
The paper is organized into three sections. In the first section, we present some preliminaries. We define what we mean by a complete conformal metric on $\mathbb{D}$ whose derivatives of all orders are uniformly bounded with respect to the hyperbolic metric. We also recall the definition of quasi-symmetric maps and some useful properties of them. The second section concerns the hyperbolic space. We begin with the definition and basic properties of quasi-circles, then define quasi-Fuchsian manifolds, since the proof of Theorem \ref{first} relies on approximation arguments using lifts of embedded surfaces in such manifolds. We then discuss the geometry of immersed surfaces in hyperbolic space and conclude with the proof of the main theorem. The third section concerns anti-de Sitter manifolds. We start by defining the space and presenting some of its models, then introduce globally hyperbolic maximal compact manifolds, which are Lorentzian analogues of quasi-Fuchsian manifolds, and finally prove Theorem \ref{secend}.\\
The proofs of Theorem \ref{first} and Theorem \ref{secend} follow the same general approach, although the techniques differ due to the change of ambient space.
\section{acknowledgment}
 I would like to express my gratitude to Jean-Marc Schlenker for the valuable discussions and insightful comments during my last visit to Luxembourg.
\section{preliminaries}
In this section we introduce the notions and some useful lemmas that we will use later though the paper.\\
\subsection{Conformal metrics with negative curvature on $\mathbb{D}$}
We begin by giving a precise definition of a conformal metric whose derivatives of all orders are bounded with respect to the hyperbolic metric.
\begin{definition}\label{boundedness}
We denote by $h_{-1}$ the hyperbolic metric on the unit disc $\mathbb{D}$, written in conformal form as $h_{-1} = \frac{4\,|dz|^2}{(1 - |z|^2)^2}$. Let $h = e^{2\rho} h_{-1}$ be a complete conformal metric on $\mathbb{D}$, where $\rho : \mathbb{D} \to \mathbb{R}$ is a smooth function.\\
We say that $h$ has bounded derivatives if, for every $p\in\mathbb{N}$, there exists a constant $M_p>0$ such that all derivatives of $\rho$ of order $p$, measured with respect to the hyperbolic metric $h_{-1}$, are uniformly bounded on $\mathbb{D}$ by $M_p$, independently of the point $z\in\mathbb{D}$.\\
In terms of partial derivatives, this means that for any $p \in \mathbb{N}$ and for every multi-index $\alpha = (\alpha_1,\alpha_2) \in \mathbb{N}^2$ with $|\alpha| = p$, we have :
$$
\left| D^\alpha \rho(z) \right| \leq M_p (1 - |z|^2)^p \quad \text{for all } z \in \mathbb{D},
$$
where $D^\alpha = \partial_x^{\alpha_1} \partial_y^{\alpha_2}$ and $z = x + iy$. The factor $(1 - |z|^2)^p$ reflects that the derivatives are measured using the norm induced by the hyperbolic metric.\\
\end{definition}
Note that any conformal metric on $\mathbb{D}$ that has constant negative sectional curvature satisfies this condition. Also, note that any metric which is invariant under a Fuchsian representation $\rho: \pi_{1}(S) \to \PSL(2,\mathbb{R})$, where $S$ is a closed hyperbolic surface, satisfies this boundedness condition.\\
Throughout the paper, we will need to approximate a given metric by metrics invariant under Fuchsian representations, so we state the following lemma.
\begin{lemma}\cite[Lemma 6.2]{curvature-abder}\label{curvature-abder}
 Let $a,b, r > 0$. Let $K : \mathbb{D} \to \left(-a,-b\right)$ be a smooth function such that any derivative of it at order $p$ is bounded by some $M_p > 0$ uniformly on $\mathbb{D}$. Let $\rho_{n} : \pi_{1}(S_{n}) \to \PSL(2,\mathbb{R})$ be a sequence of Fuchsian representations that have injectivity radius growing to $\infty$. Then there exists a sequence of smooth functions $K_{n} : \mathbb{D} \to \left(-a-r,-b+r \right)$, such that each $K_{n}$ is $\rho_{n}$-equivariant, $K_{n}$ converge $C^{\infty}$ on compact subsets to $K$, and each derivative of order $p$ of $K_{n}$ is bounded on the disc $\mathbb{D}$ by some $M'_p$, where $M'_p$ does not depend on $n$.   
\end{lemma}
\begin{proof}
 The proof is exactly similar to the proof of \cite[Lemma 6.2]{curvature-abder}. The reader needs only to notice that all the steps of the proof follow by  changing $\frac{-1}{\epsilon}$ to $-b$ and $-1-\epsilon$ to $-a$.    
\end{proof}
Lemma \ref{curvature-abder} concerns the curvature of a given metric. However, when the curvature of a conformal metric on $\mathbb{D}$ is negative, it is closely related to the metric itself, as shown in the next theorem and the next lemma.
\begin{theorem}\label{metric-by-curvature}\cite{kazdan1974curvature}
Let $K : \mathbb{D} \to \mathbb{R}^-$ be a $C^{\infty}$ function, then there exists a unique complete metric $h$ on $\mathbb{D}$ which is conformal to $|dz|^2$ and has curvature equal to $K$   
\end{theorem}
It also follows that, 
\begin{lemma}\label{nathanial}\cite[Lemma 6.4
]{curvature-abder}
Let $h_{n}$ be a sequence of complete metrics on $\mathbb{D}$, and let $h$ be also a complete metric on $\mathbb{D}$. Assume that all the metrics are conformal to $ \left|dz \right|^{2}$.\\
For each $n$ we denote by $K_n$ the curvature of $h_n$. Moreover assume that $K_{h}$ the curvature of $h$, and $K_n$ for any $n$ belong to $\left [ -a,-b \right ]$ for some $a,b > 0$.\\
If $K_{n}$ converge uniformly $C^{\infty}$ on compact subsets to $K_{h}$, then $h_{n}$ converge uniformly $C^{\infty}$ on compact subsets to $h$.\\
Moreover, if there is a sequence $(M_p)_{p \in \mathbb{N}}$ such that any derivative of $K_{n}$ of order $p$ is bounded by $M_p$, then there is a sequence of positive real numbers $M'_p$ such that any derivative of $h_n$ of order $p$ is bounded by $M'_p$.
\end{lemma}

\subsection{Quasi-symmetric maps}
In this subsection we briefly remind the definition of quasi-symmetric maps. For more details about quasi-symmetric maps and their properties we refer the reader to \cite{hubbard2016teichmuller} and \cite{fletcher2006quasiconformal}. Later, we give \cite[Proposition 9.1]{bonsante2021induced} which will be essential to us thought the paper.\\
We denote $\mathbb{RP}^1 := \mathbb{R} \cup \left\{\infty \right\}$. Let $\phi: \mathbb{RP}^{1}  \to \mathbb{RP}^{1}  $ be a strictly increasing homeomorphism that satisfies $\phi(\infty) = \infty$. We say that $\phi$ is quasi-symmetric if there exists $k > 0$ such that,
$$\forall x \in \mathbb{R}, \ \forall t \in \mathbb{R}^{*}_{+}, \ \frac{1}{k} \leq\frac{\phi(x+t)-\phi(x)}{\phi(x) - \phi(x-t)} \leq k .  $$
In this case we say that $\phi$ is $k$ is quasi-symmetric, and we call $k$ by the quasi-symmetric constant of $\phi$.\\
If $\phi$ does not fix $\infty$, then we say that $\phi$ is $k$ quasi-symmetric if there exists an element $g \in PSL(2,\mathbb{R})$ (therefore many elements) such that $g \circ \phi (\infty) = \infty$ and $g \circ \phi$ is $k$ quasi-symmetric.\\
Quasi-symmetric maps can be viewed as extensions of quasi-conformal maps, in the sense that a map $f$ is quasi-symmetric if and only if it is the boundary extension of a quasi-conformal map.\\
Quasi-symmetric maps that are equivariant under the action of Fuchsian representations are dense in the space of all quasi-symmetric maps, as stated in the following proposition. Before giving the proposition, we recall that a normalized quasi-symmetric map is a quasi-symmetric map $f$ satisfying $f(i) = i$ for $i = 0, 1, \infty$.
  
\begin{proposition}\cite[Proposition 9.1]{bonsante2021induced}\label{approx-the-quasi}
Let $f$ be a normalized quasi-symmetric map. There is a sequence of equivariant normalized uniformly quasi-symmetric maps, $\rho_{n}^{+}, \rho_{n}^{-}: \pi_{1}(S_{n}) \to \PSL(2,\mathbb{R})$, that converge to $f$. Here, $S_{n}$ is a sequence of closed surfaces with genus $g_{n}$ going to $\infty$, and $\rho_{n}^{+}, \rho_{n}^{-}$ are a sequence of Fuchsian representations whose injectivity radius go to $\infty$.
\end{proposition}

\section{Hyperbolic space}
In this section we provide a proof of Theorem \ref{first}. We begin with some preliminaries on quasi-circles, then introduce the definition of quasi-Fuchsian manifolds, since lifts of embedded surfaces in such manifolds will be used to approximate the surface to be realized in Theorem \ref{first}. After that, we discuss immersed surfaces in hyperbolic space, following the theorems of Labourie (see \cite{immersion-labourie}). Finally, before proving Theorem \ref{first}, we recall the notion of the visual metric of a convex subset in the hyperbolic space. Our approach in this section follows the same steps as in Sections 2 and Section 3 of \cite{bonsante2021induced}.
\subsection{Quasi-circles in the ideal boundary of $\mathbb{H}^3$}
We begin by giving the definition of a quasi-conformal map. For more details see for example\cite{lehto1973quasiconformal}, \cite{hubbard2016teichmuller},\cite{fletcher2006quasiconformal}.
\begin{definition}
    Let $X$ and $Y$ be Riemann surfaces (not necessarily compact). Let\\ $f: X \to Y$ be an orientation preserving diffeomorphism. We define the Beltrami differential $\mu = \mu(f)$ by the equation $\frac{\partial f}{\partial \Bar{z}} = \mu \frac{\partial f}{\partial z}$. We say that $f$ is $K$ quasi-conformal if the dilatation number $K(f) = \frac{1 + \left| \mu \right|_{\infty}}{1 - \left| \mu \right|_{\infty}}$ is less than or equal to $K$. 
\end{definition}
Note that we don't need $f$ to be a $C^{1}$ diffeomorphism to define the notion of quasi-conformal maps. In fact, all we need is for $f$ to be a homeomorphism between $X$ and $Y$ that has derivatives in the sense of distribution that are $L^{2}$.\\
Now we define quasi-circle.
\begin{definition}
A Jordan curve $C \subset \CP^1$ is called a $K$ quasi-circle, if it is the image of $\RP^1$ by a quasi-conformal map $f: \CP^1 \to \CP^1$.    
\end{definition}
We can parametrize quasi-circles by quasi-symmetric maps via the following process.\\
Let $C \subset \CP^1$ be a quasi-circle. Assume that $C$ passes through $\{0,1,\infty\}$ (which can always be achieved by applying an element of $\PSL(2,\mathbb{C})$ to $C$).  
Let $\Omega^{\pm}$ be the connected components of $\CP^1 \setminus C$. Each $\Omega^{\pm}$ is a topological disk, so by the Riemann uniformization theorem there exist holomorphic isomorphisms $\phi^{\pm}: \Omega^{\pm} \to \mathbb{D}$.  
By Carath\'eodory's theorem, each $\phi^{\pm}$ extends continuously to a homeomorphism $\partial\phi^{\pm} : C \to S^1$. Up to applying elements of $\PSL(2,\mathbb{R})$, we may assume that $\phi^{\pm}(x) = x$ for any $x \in \{0,1,\infty\}$. The map $\partial\phi^{-} \circ (\partial \phi^{+})^{-1}: \RP^1 \to \RP^1$ is quasi-symmetric.\\  
It is a classical fact that there is a one-to-one correspondence between quasi-symmetric maps (up to action of $\PSL(2,\mathbb{R})$) and quasi-circles in $\CP^1$ (up to action of $\PSL(2,\mathbb{C})$) (see for example \cite{conformal-welding}). This process is called \emph{conformal welding}.

\subsection{Quasi-Fuchsian manifolds and hyperbolic ends}\label{quasi-Fuchsian-manifolds}
We denote by $\mathbb{H}^3$ the hyperbolic space, and by $\partial_{\infty}\mathbb{H}^3$ its ideal boundary. Recall that $\partial_{\infty}\mathbb{H}^3$ can be identified with $\CP^1$. We also denote by $\mathbb{H}^2$ the hyperbolic plane, and recall that its ideal boundary $\partial_{\infty}\mathbb{H}^2$ can be identified with $\RP^1$ (also $S^1$). Since our main theorem concerns hyperbolic space, we will denote $\CP^1$ and $\RP^1$ by $\partial_{\infty}\mathbb{H}^3$ and $\partial_{\infty}\mathbb{H}^2$, respectively.\\
A quasi-Fuchsian manifold is the quotient of $\mathbb{H}^3$ by a quasi-Fuchsian group $\Gamma$, where a quasi-Fuchsian group is a torsion-free, discrete subgroup of $\PSL(2,\mathbb{C})$ whose limit set is a quasi-circle.\\
Assume that a quasi-circle $C \subset \partial_{\infty}\mathbb{H}^3$ is invariant under the action of a quasi-Fucshian group $\Gamma$, and  that $\Gamma$ is isomorphic to $\pi_1(S)$ where $S$ is a closed hyperbolic surface. By Bers simultaneous uniformization theorem (see \cite{bers}), it follows that the conformal welding map that corresponds to $C$ is equivariant under the action of two Fuchsian representations $\rho_1,\rho_2 : \pi_1(S) \to \PSL(2,\mathbb{R})$ ($f$ is $\rho_1,\rho_2$ equivariant if $\rho_1 f= f \rho_2$). By the uniqueness of the conformal welding map, it follows that reciprocally if the conformal welding map that correspondence to $C$ is equivariant by Fuchsian representations, then the quasi-circle is invariant under the action of a quasi-Fuchsian representation.\\
In particular from Proposition \ref{approx-the-quasi} we get that any quasi-circle $C$ in $\partial_{\infty}\mathbb{H}^3$ is the Hausdorff limit of quasi-circles $C_n$ such that each $C_n$ is invariant under the action of a quasi-Fuchsian representation $\rho_{n}: \pi_1(S_n) \to PSL(2,\mathbb{C)}$. Moreover, since quasi-circles are dense in the space of Jordan curves with respect to the Hausdorff topology, it follows that every Jordan curve is also the Hausdorff limit of a sequence of quasi-circles $C_n$ , where each $C_n$ is invariant under the action of a quasi-Fuchsian representation $\rho_n : \pi_1(S_n) \to PSL(2,\mathbb{C})$.\\
For a Jordan curve $C$, we denote by $CH(C)$ its convex hull, if $C$ is invariant by a quasi-Fuchsian group $\Gamma$, we denote $Q := \mathbb{H}^3 / \Gamma$, and we denote $C(Q) := CH(C) / \Gamma$, we call $C(Q)$ the convex core of $Q$.\\
The set $Q \setminus C(Q)$ has two connected components, each of these two components is called a hyperbolic end. We draw the reader attention that a hyperbolic end is a more general notion (See \cite[Definition 1.2]{Labourie-curvature}), but for us we will deal only with hyperbolic ends that comes from a quasi-Fuchsian manifold.\\
The key theorem of this subsection is the following.
\begin{theorem}\cite[Theorem 1]{Labourie-curvature}\label{Labourie-curvature}
Let $B$ be a hyperbolic end. Let $S$ be a closed surface of genus $g \geq 2$. Let $K$ be a $C^\infty$ function defined on $S$ with values in $(-1,0)$.\ 
Then there exists a unique convex incompressible surface $S$ in $B$ such that the curvature of the induced metric on $S$ from $B$ is equal to $K$. 
\end{theorem}
 Theorem \ref{Labourie-curvature} implies that if a quasi-circle $C$ is invariant under the action of a quasi-Fuchsian representation $\rho : \pi_1(S) \to \PSL(2,\mathbb{C})$, and if $K : \mathbb{D} \to (-1,0)$ is a smooth function invariant under a Fuchsian representation $\rho' : \pi_1(S) \to \PSL(2,\mathbb{R})$, then there is an isometric embedding $V : (\mathbb{D},h) \to \mathbb{H}^3$ such that the curvature of $h$ is equal to $K$ and such that the surface $\Sigma := V(\mathbb{D})$ has ideal boundary $\partial_{\infty}\Sigma = C$. 

\subsection{Immersed surfaces in the hyperbolic space}
In this section, we present some results that will be useful later to show that the surfaces we construct converge to the one we aim to realize in Theorem \ref{first}. We want to draw the reader's attention to the fact that, in the next theorem, the surface $S$ is not necessarily compact, and the metrics induced on $S$ by the immersions into $\mathbb{H}^3$ are not necessarily complete.
\begin{theorem}\cite[Theorem D]{immersion-labourie}\label{immersion-labourie}
Let $f_n : S \to \mathbb{H}^3$ be a sequence of immersions of a surface $S$ such that the pullback $f_n^{\ast}(h)$ of the hyperbolic metric $h$ converges smoothly to a metric $g_0$, and such that all the sectional curvatures of $f_n^{\ast}(h)$ are bigger than $-1+\epsilon$, for some $\epsilon > 0$. Also assume that there exists a point $x_0 \in S$ such that $f_n(x_0)$ converges and the $1$-jets of $f_n$ at $x_0$ also converge.  If the integral of the mean curvatures is uniformly bounded, then a subsequence of $f_n$ converges smoothly to an isometric immersion $f$ such that $f^{\ast}(h) = g_0$.
\end{theorem}
Also, the following lemma fit perfectly our case,
\begin{lemma}\cite[Lemma 3.7]{bonsante2021induced}\label{bounded-mean-curvature}
Let $f : S \to \mathbb{H}^3$ be a convex embedding and let $R$ be the extrinsic diameter of $h(S)$. Denote by $H$ the mean curvature and by $da$ the area form induced by $f$. Then we have
\[
\int_S H \, da \;<\; \frac{1}{\sinh1}A(R+1),
\]
where $A(r)$ denotes the area of the sphere of radius $r$ in the hyperbolic space.
    
\end{lemma}
The following lemma shows that the principal curvatures of an embedded convex surface depend only on its induced metric, and not on the surface itself, the case of surfaces of constant curvature was treated in\cite[Proposition 3.8]{bonsante2021induced}. Here we adapt the proof to the case of variable curvature.
\begin{lemma}\label{bounded-principal-curvatures}
 Let $C$ be a Jordan curve in $\partial_{\infty}\mathbb{H}^3$. Given a conformal metric $h$ on $\mathbb{D}$ that has curvature in $(-1+\epsilon,-\epsilon)$ and its hyperbolic derivatives at any order $p$ are bounded with some $M_p > 0$. Assume that there is an isometric embedded surface $V: (\mathbb{D},h) \to \mathbb{H}^3$ such that $\partial_{\infty}S = C$ (where $S = V(\mathbb{D})$). Then there exists $N > 1$ that depends only on $\epsilon$ and $(M_p)_{p \in \mathbb{N}}$ such that the principal curvatures of $S$ are in the interval $(\frac{1}{N},N)$.
\end{lemma}
\begin{proof}
 Since the curvature of $S$ belongs to the interval $(-1+\epsilon,-\epsilon)$, it follows that the product of the principal curvatures of $S$ belongs to $(\epsilon, 1-\epsilon)$. Therefore, it suffices to show that the principal curvatures of $S$ are uniformly bounded from above. We argue by contradiction. Assume there exists a sequence of such surfaces $S_n$ and a sequence of points $p_n \in S_n$ such that $\kappa_{n}(p_n)$, the largest principal curvature of $S_n$ at $p_n$, goes to $\infty$.\\
Up to applying a sequence of isometries of $\mathbb{H}^3$, we may assume that $p_n = p$, a fixed point of $\mathbb{H}^3$, and that $T_{p_n}S_n$ is a fixed tangent plane.\\
Let $V_n : (\mathbb{D},h_n) \to S_n$ be the sequence of embeddings. Up to applying isometries of $\mathbb{H}^2$, we may also assume the existence of a fixed point $x_0 \in \mathbb{H}^2$ such that $V_n(x_0) = p$. Note that these isometries of $\mathbb{H}^2$ do not necessarily preserve the metrics $h_n$, but they preserve the fact that the curvature belongs to $(-1+\epsilon,-\epsilon)$ and that the derivatives of order $m$ are bounded by the constants $M_m$ given in the statement of the lemma. For this reason, we keep the same notations $V_n$ and $h_n$ after applying such isometries, since it does not affect the proof.\\
First, observe that the metrics $h_n$ converge, up to a subsequence, smoothly and uniformly on compact subsets to a metric $h$ that satisfies the same curvature bounds and derivative estimates as the $h_n$. 
Recall that by \cite{TS-Yau} (also \cite{Troyanov}) the metrics $h_n$ are uniformly bilipschitz to the hyperbolic metric on $\mathbb{D}$ (because they have uniformly bounded negatives curvatures). Then, by Theorem \ref{immersion-labourie} and Lemma \ref{bounded-mean-curvature}, there exists a neighborhood $\mathcal{U}$ of $x_0$ (take a hyperbolic ball for example since all $h_n$ are uniformly bilipchitz to the hyperbolic metric) such that, up to extracting a subsequence, the restrictions $V_{n}\mid_{\mathcal{U}}$ converge smoothly and uniformly to an immersion 
$V : (\mathcal{U},h \mid_{\mathcal{U}}) \to \mathbb{H}^3$
with $V(x_0) = p$. This contradicts the assumption that $\kappa_n(p)$ goes to $\infty$.
 
\end{proof}
\subsection{The visual metric}
Here we will use similar technics to \cite[Section 3.3]{bonsante2021induced}.\\
Let $\mathcal{C}$ be a closed convex subset of $\Bar{\mathbb{H}}^3$ (where $\Bar{\mathbb{H}}^3 := \partial_{\infty}\mathbb{H}^3 \cup \mathbb{H}^3$). We denote by $\partial_{\infty}\mathcal{C}$ the ideal boundary of $\mathcal{C}$, that is $\mathcal{C} \cap \partial_{\infty}\mathbb{H}^3$, and we denote by $\partial\mathcal{C}$ the boundary of $\mathcal{C}$ in $\mathbb{H}^3$. It is a classical fact (see for example \cite[Chapter II.1]{canary2006fundamentals}) that there is a map, called the nearest retraction map $r_{\mathcal{C}}:\mathbb{H}^3 \to \mathcal{C}$ that associate to each point $x \in \mathbb{H}^3$ the nearest point in $\mathcal{C}$. The map $r_{\mathcal{C}}$ induces a map $r_{\mathcal{C}}: \mathbb{H}^3 \setminus \mathcal{C} \to \partial \mathcal{C}$ which is 1-Lipchitz (see \cite[Chapter II.1]{canary2006fundamentals}). Moreover the map $r_{\mathcal{C}}$ extends to a retraction from $\Bar{\mathbb{H}}^3 \to \mathcal{C}$. We extend it in the following way: let $x \in \partial_{\infty}\mathbb{H}^3 \setminus{\mathcal{C}}$, we define $r_{\mathcal{C}}(x)$ to be the intersection of the smallest horosphere with $\partial \mathcal{C}$, if $x \in \mathcal{\mathcal{C}}$ then $r_{\mathcal{C}}(x) = x$. Also note that if $\mathcal{C}_n$ converge to $\mathcal{C}$ in the Hausdorff topology of $\Bar{\mathbb{H}}^3$ then $r_{\mathcal{C}_n}$ converge uniformly to $r_{\mathcal{C}}$.\\
There is a natural metric on $\partial_{\infty}\mathbb{H}^3 \setminus \mathcal{C}$ induced by the convex subset $\mathcal{C}$ which is called the horospherical metric.\\
Before giving the definition of the horospherical metric induced by $\mathcal{C}$, and which will be denoted by $I^*_{\mathcal{C}}$, let's remind that any point $x \in \mathbb{H}^3$ induces a metric $Vis_x$ on $\partial_{\infty}\mathbb{H}^3$ called the visual metric. Also, two points $x_1,x_1 \in \mathbb{H}^3$ induces the same visual metric on the tangent space to $z \in \partial_{\infty}\mathbb{H}^3$ if and only if $x_1$ and $x_2$ belong to the same horoball centered at $z$ (by this, we mean that $Vis_{x_1}(z)=Vis_{x_2}(z)$ if and only if $x_1$ and $x_2$ lie in the same horoball centered at $z$). \\
Let $\mathbb{B}$ be the space of horospheres, and let $\pi: \mathbb{B} \to \partial_{\infty}\mathbb{H}^3$ be the map that associate to each horoball its center. A natural section of $\pi$ is $\sigma_{\mathcal{C}}: \partial_{\infty}\mathbb{H}^3 \setminus \mathcal{C} \to \mathbb{B}$, where $\sigma_{\mathcal{C}}$ associate to each point $z \in \partial_{\infty}\mathbb{H}^3 \setminus \mathcal{C}$ the horoball that is tangent to $\partial\mathcal{C}$ at $r_{\mathcal{C}}(z)$.\\
For any $z \in \partial_{\infty}\mathbb{H}^3 \setminus \mathcal{C}$, we endow $T_{z}\partial_{\infty}\mathbb{H}^3$ with the visual metric induced from the points of $\sigma_{\mathcal{C}}(z)$. By doing the same process with all the points of $\partial_{\infty}\mathbb{H}^3 \setminus \mathcal{C}$ we obtain a metric $I^{*}_{\mathcal{C}}$ on $\partial_{\infty}\mathbb{H}^3 \setminus \mathcal{C}$ which is called the horospherical metric.\\
The following remark will be useful for us later,
\begin{remark}\cite[Remark 3.9]{bonsante2021induced}\label{remark}
Note that if $\mathcal{C}_1 \subset \mathcal{C}_2$ then $I^*_{\mathcal{C}_1}\mid_{\partial_{\infty}\mathbb{H}^3\setminus \partial_{\infty}\mathcal{C}_2} \leq I^*_{\mathcal{C}_2}$. For us the most important fact is that if $\mathcal{C}_1$ and $\mathcal{C}_2$ share the same ideal boundary (that is $\partial_{\infty}\mathcal{C}_1 = \partial_{\infty}\mathcal{C}_2$) and $I^*_{\mathcal{C}_1} \leq I^*_{\mathcal{C}_2}$ then $\mathcal{C}_1 \subset \mathcal{C}_2$.    
\end{remark}

\begin{figure}
    \centering
    \includegraphics[width=0.5\linewidth]{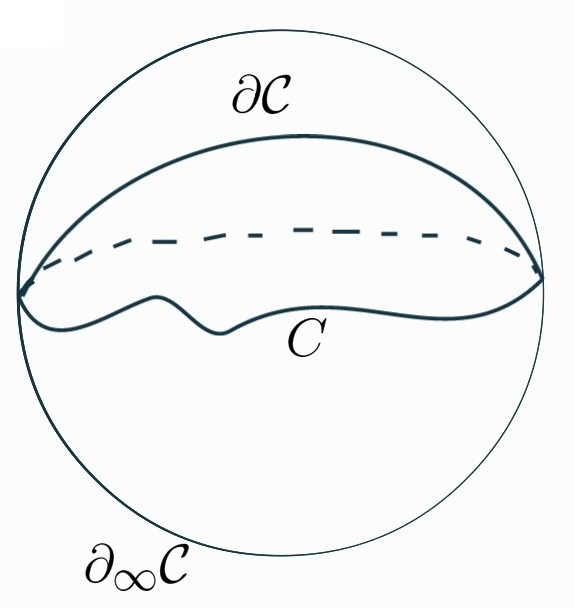}
    \caption{For us, the convex set $\mathcal{C}$ will always be a convex subset of $\overline{\mathbb{H}}^3$ whose boundary in $\mathbb{H}^3$ is a smooth surface $\widetilde{S}$, with ideal boundary $\partial_{\infty}\widetilde{S}$ equal to a Jordan curve $C$. Moreover, $\partial_{\infty}\mathcal{C}$ is given by $C \cup \Omega^{-}$, where $\Omega^{-}$ is the component of $\partial_{\infty}\mathbb{H}^3 \setminus C$ lying on the convex side of $\widetilde{S}$.}
    \label{Figure-convex-C}
\end{figure}

It is also worth noticing that in the case when $\partial \mathcal{C} = \partial^+CH(C)$ (that is the upper boundary component of the convex hull of a Jordan curve $C$) then $I^*_{\mathcal{C}}$ is equal to the Thurston metric on $\partial_{\infty}\mathbb{H}^3 \setminus \partial_{\infty}C$ (for more details see \cite{Canary-Bridgman}).\\
For simplicity, we denote $\partial_{\infty}\mathbb{H}^3 \setminus \partial_{\infty}\mathcal{C} := \Omega^{+}_{\mathcal{C}}$, 
and thus $\partial_{\infty}\mathcal{C} = \Omega^{-}_{\mathcal{C}}$. 
This notation will be useful since we will consider a convex surface $\tilde{S}$ spanning a Jordan curve $C$ at infinity. In this case, $\partial_{\infty}\mathbb{H}^3 \setminus C$ consists of two connected components, denoted by $\Omega^+$ and $\Omega^-$. Unless stated otherwise, $\mathcal{C}$ will refer to the convex subset bounded by $\Omega^-$ and $\tilde{S}$ (see Figure \ref{Figure-convex-C}). Note also that in this setting, $\partial_{\infty}\mathbb{H}^3 \setminus \mathcal{C} = \Omega^+$.
\subsection{The proof of the main theorem}
The next lemma is a key step in proving Theorem \ref{first}. It ensures that the surface $\widetilde{S}$ can be approximated by a sequence $\widetilde{S}_n$, which we will later show converges to $\widetilde{S}$.
\begin{lemma}\label{approximate-hyperbolic-surfaces}
Let $C$ be a Jordan curve in $\partial_{\infty}\mathbb{H}^3$, and let $h$ be a complete, conformal metric on $\mathbb{D}$ whose curvature belongs to $(-1+\epsilon,-\epsilon)$ and whose hyperbolic derivatives are bounded at every order.\\
Then there exists a sequence of quasi-circles $C_n \subset \partial_{\infty}\mathbb{H}^3$ converging to $C$ in the Hausdorff topology, and a sequence of isometric embeddings
$$
V_n : (\mathbb{D},h_n) \to \tilde{S}_n \subset \mathbb{H}^3,
$$
such that the boundary at infinity of $\tilde{S}_n$ satisfies $\partial_{\infty}\tilde{S}_n = C_n$. Moreover, the metrics $h_n$ have curvature in $(-1+\epsilon,-\epsilon)$, their hyperbolic derivatives are uniformly bounded at every order, and $h_n$ converge in $C^{\infty}$ on compact subsets to $h$.
\end{lemma}
\begin{proof}
By Proposition \ref{approx-the-quasi} there exists a sequence $C_n$ of quasi-circles converging in the Hausdorff topology to the Jordan curve $C$. Moreover, each $C_n$ is invariant under a quasi-Fuchsian representation $\rho_n : \pi_1(S_n) \to \PSL(2,\mathbb{C})$.\\ 
Let $M_n$ be the quasi-Fuchsian manifold $\mathbb{H}^3/\rho_n$, which is homeomorphic to $S_n \times \mathbb{R}$. We denote by $M_n^+$ the upper connected component of $M_n \setminus C(M_n)$, where $C(M_n)$ is the convex core of $M_n$.\\
Let $K_h$ be the curvature of $h$. By Lemma \ref{curvature-abder}, there exists a sequence of functions $\bar{K}_n : S_n \to (-1+\epsilon,-\epsilon)$ whose lifts $K_n : \mathbb{D} \to (-1+\epsilon,-\epsilon)$ converge to $K_h$ in $C^{\infty}$. Moreover, all the functions $K_n$ have uniformly bounded derivatives at any order with respect to the hyperbolic metric.\\
By Theorem \ref{Labourie-curvature}, there exists an embedded surface $S_n \subset M_n^+$ whose curvature is $\bar{K}_n$. The lift of the surface $S_n$ gives to an embedding $V_n : \mathbb{D} \to \mathbb{H}^3$ 
such that the curvature of the induced metric of $\tilde{S}_n := V_n(\mathbb{D})$ is $K_n$. By Theorem \ref{metric-by-curvature}, there exists a unique complete conformal metric $h_n$ on $\mathbb{D}$ with curvature $K_n$, and by Lemma~\ref{nathanial} the sequence $h_n$ converges to $h$ smoothly uniformly on compact subsets (because by hypothesis $K_n$ converge to $K_h$). Also by Lemma \ref{nathanial}, the metrics $h_n$ have uniformly bounded hyperbolic derivatives at every order.    
\end{proof}
The next step is to show that the surfaces $\tilde{S}_n$ converge to a surface $\tilde{S}$ which is isometric to $(\mathbb{D},h)$ and satisfies $\partial_{\infty}S = C$.\\
We will need the following lemma. Before giving the lemma, If $\partial \mathcal{C}$ is of regularity $C^2$, let's denote by $I_{\mathcal{C}}$ the first fundamental form on $\mathcal{C}$ (that is the induced metric), we denote by $\II_{\mathcal{C}}$ the second fundamental form on $\partial \mathcal{C}$, and by $\III_{\mathcal{C}}$ the third fundamental form on $\partial \mathcal{C}$. 
\begin{lemma}\cite{SH02}\label{SH02}
  \begin{enumerate}
    \item If $\mathcal{C}_s$ is the set of points at distance less than or equal to $s$ from $\mathcal{C}$, then $\mathcal{C}_s$ is a convex set and $I^*_{\mathcal{C}_{s}} = e^sI^*_{\mathcal{C}}$.

    \item If $\partial\mathcal{C}$ is of class $C^2$, then $r_{\mathcal{C}} : \op \to \partial\mathcal{C}$ is a $C^1$-diffeomorphism and $$(r_{\mathcal{C}}^{-1})^*(I^*_{\mathcal{C}}) = I_{\mathcal{C}} + 2\II_{\mathcal{C}} + \III_{\mathcal{C}}.$$

    \item If $\partial \mathcal{C}$ is smooth, then the curvature of $I^*_{\mathcal{C}}$ at $z \in \op$ is 
    $$K(z) = \dfrac{K(r_{\mathcal{C}}(z))}{(1 + \kappa_1(r_{\mathcal{C}}(z)))(1 + \kappa_2(r_{\mathcal{C}}(z)))},$$  
    where $K$ is the intrinsic curvature of $(\partial \mathcal{C}, I_{\mathcal{C}})$ and $\kappa_1, \kappa_2$ denote the principal curvatures.
\end{enumerate}
 
\end{lemma}
As a consequence, we obtain the following proposition. The case of constant Gaussian curvature was already proved in \cite[Proposition 3.11 and Lemma 3.13]{bonsante2021induced}.
\begin{proposition}\label{bilipchitz-estimate}
 Let $h_{\op}$ be the hyperbolic metric on $\op$ compatible with its conformal structure. Assume that $(\partial\mathcal{C},I_{\mathcal{C}})$ is isometric to $(\mathbb{D},h)$ where $h$ has curvature in $(-1+\epsilon,-\epsilon)$ and the derivatives of $h$ are bounded at any order. Then:
 \begin{itemize}
     \item $h_{\op}$ is billipchitz to $I^*_{\mathcal{C}}$.
     \item The map $r_{\op}: (\op, I^*_{\op}) \to (\partial\mathcal{C},I_{\mathcal{C}})$ is bilipchitz.
     \item The map $r_{\op}: (\op, h_{\op}) \to (\partial\mathcal{C},I_{\mathcal{C}})$ is bilipchitz.
 \end{itemize}
  Moreover, the billipchitz constants depend only on the curvature of $h$ and the bounds of its derivatives.    
\end{proposition}
\begin{proof}
To prove the first point, first, note that by Lemma \ref{bounded-principal-curvatures} the principal curvatures of $\partial\mathcal{C}$ are bounded with bounds that depend only on the curvature of $h$ and the bounds of its derivatives. Then, Lemma \ref{SH02}(3) implies the existence of $M > 0$, that also depends only on the curvature of $h$ and the bounds of its derivatives, such that the curvature of $I^*_{\op}$ is bounded by $-M$ and $-\frac{1}{M}$. Also by Lemma \ref{SH02}(2) the metric $ I^*_{\op}$ is complete. It follows from a Theorem of Yau \cite{TS-Yau} (see also the first theorem of \cite{Troyanov}) that the metrics $I^*_{\op}$ and $h_{\op}$ are $M$ bilipchitz.\\
To prove the second point, if suffices to use Lemma \ref{SH02} (2) and Lemma \ref{bounded-principal-curvatures}. Indeed, Lemma \ref{bounded-principal-curvatures} implies that the principal curvatures of $\partial \mathcal{C}$ are bounded, with bounds that depend only on the curvature of $h$ and the bounds of its derivatives. Also,  Lemma \ref{SH02} (2) implies that the pull back of $I_{\partial\mathcal{C}}$ by $r_{\op}$ is a linear combination of the first, the second, and the third fundamental forms of $\partial\mathcal{C}$, so the statement follows.\\
The third point is a direct consequence of the first and the second.

\end{proof}

\begin{corollary}\label{extend-hyperbolic}
The isometric embedding $V : (\mathbb{D},h) \to \mathbb{H}^3$ extends continuously to $\partial V: S^1 \to C$    
\end{corollary}
\begin{proof}
Let $\phi_{\op}: \op \to \mathbb{H}^2$ be the uniformization. By Carath\'eodory Theorem, it extends to a homeomorphism $\partial \phi_{\op}: C \to S^1$. Also note that the map $r_{\op}:\op \to \partial\mathcal{C}$ extends continuously to the identity $\mathrm{Id}: C \to C$.\\
By Proposition \ref{bilipchitz-estimate}, the map $\phi_{\op}\circ r^{-1}_{\op} \circ V : (\mathbb{D},h) \to \mathbb{H}^2$ is bilipschitz, so it extends to a homeomorphism $S^1 \to \partial_{\infty}\mathbb{H}^2$. Since $\phi_{\op}$ extends to a homeomorphism $\partial \phi_{\op}: C \to \partial_{\infty}\mathbb{H}^2$ and $r_{\op}$ extends to the identity $C \to C$, it follows that $V$ extends to a homeomorphism $\partial V: S^1 \to C$.
\end{proof}
The following lemma was proved in the constant Gaussian curvature case in \cite[Lemma 4.5]{bonsante2021induced}. Here, we extend it to a broader class of metrics with variable curvature.
\begin{lemma}\label{distance}
 There is a constant $R$ that depends only on the curvature of $h$ and the bounds of its derivatives such that if $(\partial\mathcal{C}, I_{\mathcal{C}})$ is isometric to $(\mathbb{D},h)$ then $\partial\mathcal{C}$is in an $R$ neighborhood of $\partial^{+}\CH(C)$.   
\end{lemma}
\begin{proof}
 By Proposition \ref{bilipchitz-estimate} the identity map $id: (\op, I^*_{\mathcal{C}}) \to (\op,h_{\op}) $ is a bliptchitz map, and the bilipchitz constant depends only on the curvature of $h$ and the bounds of its derivatives. Denote by $I^*_{Th}$ the Thurston metric on $\op$ (that is the horospherical metric on $\op$ which correspondence to $\partial^+CH(C)$), Bt \cite[Lemma 3.1]{thurston-estimate} we have that $h_{\op} \leq 2I^*_{Th}$. This implies the existence of $R > 0$ that depends only on the curvature of $h$ and the bounds of its derivatives such that $I^*_{\op} \leq e^RI^{*}_{Th}$. By the first point of Proposition \ref{SH02} and by Remark \ref{remark} we conclude that $\partial\mathcal{C}$ is at distance at most $R$ from $\partial^+CH(C)$.    
\end{proof}
Then under the hypothesis on $h_n$, $h$, $C_n$ and $C$ given in this section, the following  lemmas hold.
\begin{lemma}
There exists $x_0 \in \mathbb{D}$ such that $V_n(x_0)$  is in a compact subset of $\mathbb{H}^3$ (up to extracting a subsequence). 
\end{lemma}
\begin{proof}
Consider the normalized map $\Phi_n = V_n^{-1} \circ r_{\Omega^+_{C_n}} \circ U_n : \mathbb{D} \to \mathbb{D}$.  
If we equip $\mathbb{D}$ with the conformal hyperbolic metric, by Proposition \ref{bilipchitz-estimate} the maps $\Phi_n$ are uniformly bilipschitz.\\ 
By \cite[Lemma 4.8]{bonsante2021induced}, this implies that there exists $x_0 \in \mathbb{D}$ such that $\Phi_n(x_0)$ remains in a compact subset of $\mathbb{H}^2$.\\  
Then $V_n(\Phi_n(x_0)) = r_{\Omega^+_{C_n}} \circ U_n(x_0)$ belongs to a compact subset of $\mathbb{H}^3$, since $U_n \to U$, $r_{\Omega^+_{C_n}} \to r_{\Omega^+_{C}}$, and $\partial\mathcal{C}$ lies at a bounded distance from $\partial^+CH(C)$.\\  
This concludes the proof, as $\Phi_n(x_0)$ stays in a compact subset of $\mathbb{D}$ and thus, up to extracting a subsequence, converges to some point in $\mathbb{D}$.
\end{proof}

\begin{lemma}\cite[Lemma 4.6]{bonsante2021induced}\label{hyperbolic-immersions-conerge-to-an-embedding}
There exists $x_0 \in \mathbb{D}$ such that $V_n(x_0)$  is in a compact subset of $\mathbb{H}^3$ (up to extracting a subsequence). Then, the isometric embeddings $V_n:(\mathbb{D},h_n) \to \mathbb{H}^3$ converge (up to extracting a subsequence) to an isometric embedding $V:(\mathbb{D},h) \to \mathbb{H}^3$ that extends continuously to a homeomorphism $\partial V : S^1 \to C$.   
\end{lemma}
\begin{proof}
Note that, by hypothesis, the quasi-circles $C_n$ converge in the Hausdorff topology of $\overline{\mathbb{H}}^3$ to the Jordan curve $C$.\\  
The convex domains $\mathcal{C}_n$ constructed in Lemma \ref{approximate-hyperbolic-surfaces} converge, up to extracting a subsequence, to a convex subset $\mathcal{C}$. Since $C_n = \partial_{\infty}\mathcal{C}_n$ converge to $\partial_{\infty}\mathcal{C}$, it follows that $\partial_{\infty}\mathcal{C} = C$. Moreover, $\partial^+CH(C_n)$ converge to $\partial^+CH(C)$.\\  
We also denote by $\Omega^+_{C_n}$ and $\Omega^+_{C}$ the connected components of $\partial_{\infty}\mathbb{H}^3 \setminus C_n$ and $\partial_{\infty}\mathbb{H}^3 \setminus C$, respectively, lying on the concave side of $\partial^+CH(C_n)$ and $\partial^+CH(C)$.  
Let $U_n : \mathbb{D} \to \Omega^+_{C_n}$ and $U: \mathbb{D} \to \Omega^+_{C}$ be the Riemann uniformization maps.\\  
Since $C_n$ converge to $C$, we may assume, without loss of generality, the existence of three points (for example $0,1,\infty$) through which both $C_n$ and $C$ pass. Indeed, choose three distinct points $x,y,z$ in $C$ and sequences of points $x_n,y_n,z_n$ in $C_n$ such that $x_n \to x$, $y_n \to y$, and $z_n \to z$. Then apply isometries $A_n$ on $C_n$ such that $A_n(x_n) = x$, $A_n(y_n) = y$, and $A_n(z_n) = z$. The isometries $A_n$ converge to the identity since $x_n,y_n,z_n \to x,y,z$.\\  
From now on, we assume that all the isometries and uniformization maps are normalized, in the sense that each of the points $i=0,1,\infty$ is fixed by these maps.\\  
Since each $\partial \mathcal{C}_n$ is at distance at most $R$ from $\partial^+CH(C_n)$, it follows that $\partial \mathcal{C}$ is at distance at most $R$ from $\partial^+CH(C)$.  

It also follows that $\partial \mathcal{C}$ is a topological disc. Indeed, since $\partial \mathcal{C}$ lies within distance $R$ of $CH(C)$, it cannot intersect $\Omega^{+}_{CH(C)}$. Therefore, the projection map $r_{\Omega^{+}_{C}}$ is a homeomorphism on $\partial \mathcal{C}$.
 As a consequence, $\Omega^+_{\mathcal{C}} = \Omega^+_{CH(C)}$.  
By Theorem \ref{immersion-labourie} and Lemma \ref{bounded-mean-curvature}, the isometric embeddings $V_n : (\mathbb{D},h_n) \longrightarrow \partial \mathcal{C}_n$  
converge smoothly, up to subsequence, on compact subsets to an immersion  
$V : (\mathbb{D},h) \longrightarrow \partial \mathcal{C}$.\\  
Since $(\mathbb{D},h)$ is complete, the map $V$ is a covering. As $\partial \mathcal{C}$ is simply connected, $V$ is a homeomorphism. We deduce that $V$ is an isometric embedding of $(\mathbb{D},h)$ into $\mathbb{H}^3$.\\  
Finally, by Lemma \ref{extend-hyperbolic}, the map $V$ extends continuously to a homeomorphism  
$\partial V : S^1 \longrightarrow C$.  
\end{proof}
\section{Anti-de sitter space}
In this section we give a proof of Theorem \ref{secend}. We will use the same approach as in \cite[Section 6 and Section 7]{bonsante2021induced} and in \cite{curvature-abder}.\\
We begin by briefly recalling the definition of $\AdS^3$, mainly to fix notations.
\subsection{Preliminaries on the anti-de sitter space}\label{anti-de sitter}
The anti-de Sitter space $\AdS^3$ is a Lorentzian space, regarded as the Lorentzian analogue of the hyperbolic space $\mathbb{H}^3$. It carries a Lorentzian metric of constant sectional curvature $-1$, and any three-dimensional Lorentzian manifold with constant sectional curvature $-1$ is locally modeled on $\AdS^3$.\\
In this section, we introduce the models of the three-dimensional anti-de Sitter space $\AdS^{3}$ that are mostly used. For further details, we refer the reader to \cite{bonsante2020anti} and \cite{zbMATH05200424}.
\subsection{Hyperbloid model}
Let $q_{2,2}$ be the quadratic form defined on $\mathbb{R}^{4}$ by the formula :
$$q_{2,2}(x_1,x_2,x_3,x_4) = x^2_1 + x^2_2 - x^2_3 - x^2_4.$$
We denote by $\mathbb{R}^{2,2}$ the space $(\mathbb{R}^{4},q_{2,2})$.\\
We define the space $\mathbb{H}^{2,1}$ to be :
$$\mathbb{H}^{2,1} :=  \left\{x \in \mathbb{R}^{4}, q_{2,2}(x) = -1 \right\}.$$
Then $q_{2,2}$ induce a scalar product on each tangent space of $\mathbb{H}^{2,1}$ that has signature $(2,1)$, this makes $\mathbb{H}^{2,1}$ a Lorentzian manifold. We refer for example to \cite{bonsante2020anti} to see why $\mathbb{H}^{2,1}$ has a constant sectional curvature equal to $-1$.
\subsection{The projective model of $\AdS^{3}$}
We introduce $\AdS^{3}$, the projective model of the anti-de Sitter space, to be :
$$\AdS^{3} := \mathbb{H}^{2,1} /  \left\{\pm \right\},$$
or equivalently :
$$\AdS^{3} := \left\{\left [ x \right ] \in \mathbb{RP}^3, q_{2,2}(x) < 0  \right\}.$$
The projective model allows us to visualise better $\partial_{\infty} \AdS^{3}$, the ideal boundary of the anti-de Sitter space.
$$\partial_{\infty} \AdS^{3} := \left\{\left [ x \right ] \in \mathbb{RP}^3, q_{2,2}(x) = 0  \right\}$$
The anti-de Sitter space $\AdS^{3}$ induces a Lorentzian conformal structure on its ideal boundary $\partial_{\infty}\AdS^{3}$ (see \cite[Section 2.2]{bonsante2020anti}).
\subsection{The Lie group model}
Let $\mathcal{M}_{2,2}(\mathbb{R})$ be the space of $2 \times 2$ matrices with real coefficient. The space $(\mathcal{M}_{2,2}(\mathbb{R}), -det)$ is isometric to $\mathbb{R}^{2,2}$ via the map
\begin{align*}
 \mathbb{R}^{4} & \to \mathcal{M}_{2,2}(\mathbb{R})\\
 (x_1,x_2,x_3,x_4) & \mapsto \begin{pmatrix}
x_1 - x_3 & x_4 - x_2  \\
x_2 + x_4 & x_1 + x_3  \\
\end{pmatrix}     
\end{align*}

and under this isomorphism $\mathbb{H}^{2,1}$ is identified with the Lie group $SL(2,\mathbb{R})$ (see for example \cite[Section 2.1]{bonsante2021induced}). Under this identification it yields that the projective model $\AdS^{3}$ is identified with $\PSL(2,\mathbb{R})$ and $\partial_{\infty}\AdS^{3}$ is identified with $\left\{\left [ M \right ] \in \PSL_{2}(\mathbb{R}), det(M) = 0  \right\}$.\\
There is an explicite identification between $\partial_{\infty}\AdS^{3}$ and $\mathbb{RP}^{1} \times \mathbb{RP}^1$ via the following map :
\begin{align*}
 \partial_{\infty} \AdS^{3} & \to \mathbb{RP}^1 \times \mathbb{RP}^1\\
 \left [ M \right ] &  \to (Im(M),Ker(M)).
\end{align*}
Note that group $\PSL(2,\mathbb{R}) \times \PSL(2,\mathbb{R})$ acts on $\PSL(2,\mathbb{R})$ by left and right composition, that is :
$$(A,B).X = AXB^{-1}.$$
We refer the reader to \cite[Section 3.1]{bonsante2020anti} to see why we can identify the isometry group of $\AdS^{3}$ that preserve orientation and time orientation with $\PSL(2,\mathbb{R}) \times \PSL(2,\mathbb{R})$.\\
In general, in any Lorentzian manifold $(M,q)$, we say that a vector $v \in T_{p}M$ is :
\begin{itemize}
    \item Space-like if $q(v,v) > 0$.
    \item Light-like if $q(v,v) = 0$.
    \item Time-like if $q(v,v) < 0$.
\end{itemize}
The geodesics in $\AdS^{3}$ are obtained by the intersection of planes of $\mathbb{R}^{2,2}$ that go thought the origin with $\mathbb{H}^{2,1}$. We say that a geodesic $\alpha$ is :
\begin{itemize}
    \item Space-like if $q_{2,2}(\Dot{\alpha}) > 0$.
    \item Light-like if $q_{2,2}(\Dot{\alpha}) = 0$.
    \item Time-like if $q_{2,2}(\Dot{\alpha}) < 0$.
\end{itemize} 
We refer to \cite[Section 2.3]{bonsante2020anti} to see that $\AdS^{3}$ is time oriented.\\
Let $\Sigma$ be a smoothly immersed surface in $\AdS^3$, we say that $\Sigma$ is space-like if $\AdS^3$ induces on it a Riemannian metric.

\subsection{Globally hyperbolic maximal compact $\AdS^3$ manifolds}\label{Globally-hyperbolic}
Let $C$ be a continuous curve in $\partial_{\infty} \AdS^{3}$.
We say that $C$ is achronal (resp acausal), if for any point $p \in C$, there is a neighborhood $U$ of $p \in \partial_{\infty} \AdS^{3}$, such that $U \cap C$ is contained in the complement of the regions of $U$ which are connected to $p$ by timelike curves (resp timelike and lightlike curves).\\
We have seen that $\partial_{\infty}\AdS^{3}$ is identified with $\mathbb{RP}^1 \times \mathbb{RP}^1$. Then the graph of any homeomorphism $f: \mathbb{RP}^1 \to  \mathbb{RP}^1$ defines a curve on $\partial_{\infty}\AdS^{3}$.\\
By the work of Mess \cite{zbMATH05200424} the following definition holds.
\begin{definition}\label{quasi-circle}
An acausal curve $C \subset \partial_{\infty}\AdS^{3}$ is a quasi-circle if it is the graph of a quasi-symmetric map.    
\end{definition}\label{domain-of-dependance}
We will define the domain of dependence of a quasi-circle $C$. We say that a curve in $\AdS^{3}$ is causal if its tangent vector at any point is time-like or light-like.
\begin{definition}
Let $C \subset \partial_{\infty}\AdS^{3}$ be a quasi-circle. We define $D(C)$, the domain of dependence of $C$, to be:
$$
D(C) := \left\{\, p \in \AdS^{3} \;\middle|\; p \text{ is connected to } C \text{ by no causal path} \,\right\}.
$$

Equivalently, $D(C)$ (see \cite{BB09}) is the unique maximal (in the sense of inclusion) open convex subset whose boundary at infinity is equal to $C$. 
\end{definition}
For any quasi-circle $C$, the domain of dependence $D(C)$ is contained in a unique affine chart. Moreover, it is the maximal convex subset of $\AdS^{3}$ whose topological closure in $\AdS^{3} \cup \partial_{\infty}\AdS^{3}$ contains $C$.\\
For us, a globally hyperbolic maximal compact $\AdS^3$ manifold will be the quotient of some $D(C)$ by a representation $\rho : \pi_1(S) \to \PSL(2,\mathbb{R}) \times \PSL(2,\mathbb{R})$ of the form $\rho = (\rho_1,\rho_2)$ where each of $\rho_1$ and $\rho_2$ is a Fuchsian representation. We will consider only the case of globally hyperbolic maximal compact manifolds where $S$ is a closed, orientable surface of genus greater than or equal to 2.\\
An anti-de Sitter end is a connected component of $M \setminus C(M)$, where $M := D(C)/(\rho_1,\rho_2)$ is a globally hyperbolic maximal compact manifold and $C(M)$ is its convex core. The convex core is defined, as in the hyperbolic case, as the quotient of the convex hull of the quasi-circle $C$ by the action of $(\rho_1, \rho_2)$.\\
Note that Proposition \ref{approx-the-quasi} implies that any quasi-circle is the Hausdorff limit of a sequence $C_n$ of quasi-circles such that $D(C_n)$ is the lift of a globally hyperbolic maximal compact $\AdS^3$ manifold.\\
A key theorem for us is the following, this theorem will allow us to approximate the surface that we want to realize.
\begin{theorem}\cite[Theorem 2.4]{curvature-zeghib}\label{curvature-zeghib}
Let $S$ be a closed hyperbolic surface, and let $M$ be an anti-de Sitter end homeomorphic to $S \times \mathbb{R}$, assume that it is future complete. Let $K : S \to (-\infty,-1)$ be a smooth function, then there is an embedding of $S$ in $M$ such that the resulting surface is space-like and its induced metric has curvature equal to $K$.  
\end{theorem}

\subsection{Proof of the main theorem}
As in the hyperbolic case, we will use Theorem \ref{curvature-abder} to approximate the surface we want to realize. Later, our main goal will be to show that these surfaces converge in the desired sense.
\begin{proposition}\label{approximate-ads}
Let $C$ be a quasi-circle in $\partial_{\infty}\AdS^3$, and let $h$ be a complete, conformal metric on the disc $\mathbb{D}$ whose curvature $K_h$ lies in an interval of the form $(-\frac{1}{\epsilon}, -1 - \epsilon)$, with all its derivatives bounded with respect to the hyperbolic metric. Then there exists a sequence of quasi-circles $C_n \subset \partial_{\infty}\AdS^3$ converging to $C$ in the Hausdorff topology, and a sequence of complete conformal metrics $h_n$ on $\mathbb{D}$ whose curvatures lie in $(-\frac{1}{\epsilon}, -1 - \epsilon)$ and whose derivatives of all orders are uniformly bounded with respect to the hyperbolic metric, such that $h_n$ converges smoothly and uniformly on compact subsets to $h$. Moreover, for each $n$, there exists an embedding
$V_n : (\mathbb{D}, h_n) \to \AdS^3$
such that $\partial_{\infty}\Sigma_n = C_n$, where $\Sigma_n = V_n(\mathbb{D})$.   
\end{proposition}
\begin{proof}
 By Proposition \ref{approx-the-quasi} there is a sequence of quasi-circles $C_n$ invariant under the action $(\rho_n,\rho'_n) : \pi_1(S_n) \to \PSL(2,\mathbb{R}) \times \PSL(2,\mathbb{R})$ where $S_n$ is a sequence of closed hyperbolic surfaces and $\rho_n,\rho'_n$ are Fuchsian representations. Let $K_h$ be the curvature of $h$. By Lemma \ref{curvature-abder} we can construct a sequence of function $K_n : \mathbb{D} \to (-\frac{1}{\epsilon},-1-\epsilon)$ that have uniformly bounded derivatives with respect to the hyperbolic metric at any order and for each $n$, $K_n$ is invariant under $\rho_n$.\\
 This will induce a function $\Bar{K}_n : S_n \to (-\infty,-1)$. By Theorem \ref{curvature-zeghib} there is an embedded space-like surface $S_n$ in the future of $D(C_n)/(\rho_n,\rho'_n)(\pi_1(S_n))$ that has curvature $\Bar{K}_n$. This surface will lift to a space-like embedded disk $\Sigma_n$ that its induced metric has curvature $K_n$ and span $C_n$ at the boundary at infinity. Recall that there is exactly one conformal metric on $\mathbb{D}$ that has curvature $K_n$. Since $K_n$ converge to $K_h$, Lemma \ref{nathanial} implies that the induced metrics $h_n$ on $\Sigma_n$ converge smoothly uniformly on compact subsets to $h$ and all have uniformly bounded hyperbolic derivatives at any order.
\end{proof}
We recall the following lemma,
\begin{lemma}\cite[Lemma 5.6]{curvature-abder}\label{Lemma-5.6-curvature-abder}
Let $\Omega$ be a globally hyperbolic maximal compact convex subset spanning a $k$-quasicircle at infinity. Assume that the induced metrics on $\partial^\pm \Omega$ have curvatures in$(-\frac{1}{\varepsilon},\, -1-\varepsilon)$.\\
Let $V^\pm : (\mathbb{D},h^\pm) \to \partial^\pm \Omega$ be isometries, and assume that every derivative of $h^\pm$ of order $p$ is bounded by $M_p$ on $\mathbb{D}$.\\
Let $\Pi^\pm_l$ and $\Pi^\pm_r$ denote the left and right projections. Then there exists a constant $A>0$, depending only on $\varepsilon$, $k$, and $(M_p)_{p\in\mathbb{N}}$, such that the maps $\Pi^\pm_l \circ V^\pm$ and $\Pi^\pm_r \circ V^\pm$ are $A$-quasi-isometries.
\end{lemma}
The following proposition was proved in \cite{curvature-abder}, and earlier in \cite{bonsante2021induced} for the case of constant curvature. We would like to draw the reader's attention to the fact that, unlike Lemma \ref{bounded-principal-curvatures} in the hyperbolic case, here the bounds also depend on the curve $C$. Therefore, we cannot conclude Theorem \ref{secend} when $C$ is not the graph of a quasi-symmetric map.
\begin{proposition}\cite[Proposition 5.5]{curvature-abder}\label{conclude-ads}
Let $C$ be a quasi-circle, and let $h$ be a complete conformal metric on $\mathbb{D}$ whose curvature lies in the interval $\left(-\tfrac{1}{\epsilon}, -1-\epsilon\right)$, with all derivatives of any order uniformly bounded. Let $\Sigma$ be an isometrically embedded surface $V : (\mathbb{D}, h) \to \Sigma$ such that the ideal boundary $\partial_{\infty}\Sigma$ coincides with $C$.
Then the principal curvatures of $\Sigma$ lie in an interval $\left(\tfrac{1}{D}, D\right)$, where $D$ depends only on the quasi-symmetric constant of $C$, the curvature of $h$, and the bounds on the derivatives of $h$.
\end{proposition}
Next we give the following proposition from \cite{bonsante2021induced} and \cite{curvature-abder}.
\begin{proposition}\cite[Proposition 5.15]{curvature-abder}\label{last-conclusion}
Let $V_n:(\mathbb{D},h_n) \to \AdS^3$ be a sequence of smooth isometric embeddings that extend continuously to homeomorphisms $\partial V_n : S^1 \to C_n \subset \partial_{\infty}\AdS^3$. Assume that $C_n$ converge in the Hausdorff topology to a quasi-circle $C$. If all the metrics $h_n$ have curvature in an interval $(-\frac{1}{\epsilon},-1-\epsilon)$ and have uniformly bounded derivatives with respect to the hyperbolic metric at any order. Assume that $h_n$ converge smoothly uniformly on compact subsets to a complete conformal metric $h$, that also will have curvature $K_h$ that has values in the interval $(-\frac{1}{\epsilon},-1-\epsilon)$ and all derivatives bounded at any order with respect to the hyperbolic metric.\\
Then $V_n$ will converge smoothly uniformly on compact subsets to an isometric embedding $V : (\mathbb{D},h) \to \AdS^3$ that extends continuously to a homeomorphism $\partial V : S^1 \to C$   
\end{proposition}
Then we conclude the proof of Theorem \ref{secend}
\begin{proof}
Let $C \subset \partial_{\infty}\AdS^3$ be a quasi-circle. Proposition \ref{approximate-ads} ensures the existence of quasi-circles $C_n$ converging to $C$ in the Hausdorff topology, and for each $n$ there exists an isometric embedding $V_n: (\mathbb{D}, h_n) \to \Sigma_n$, such that $\partial_{\infty}\Sigma_n = C_n$. By Lemma \ref{Lemma-5.6-curvature-abder}, and by the paragraph after \cite[Lemma 5.6]{curvature-abder}, the map $V_n$ extends continuously to a homeomorphism $\partial V_n: S^1 \to C_n$.
Finally, by Proposition \ref{last-conclusion}, the isometric embeddings $V_n$ converge to an isometric embedding $V: (\mathbb{D}, h) \to \Sigma$ that extends continuously to a homeomorphism $\partial V: S^1 \to C$.  
\end{proof}

\bibliographystyle{alpha}
	\bibliography{biblo.bib}
\nocite{*}

\end{document}